# Sequential nonparametrics and semiparametrics: Theory, implementation and applications to clinical trials

## Tze Leung Lai[1] and Zheng Su[2]

*Stanford University and Genentech*

**Abstract:** One of Pranab K. Sen's major research areas is sequential nonparametrics and semiparametrics and their applications to clinical trials, to which he has made many important contributions. Herein we review a number of these contributions and related developments. We also describe some recent work on nonparametric and semiparametric inference and the associated computational methods in time-sequential clinical trials with survival endpoints.

## 1. Introduction

Sequential nonparametrics began in the 1960s with the work of Wilcoxon, Rhodes and Bradley [68] on extending Wald's [66] sequential probability ratio test (SPRT) to construct two-sample grouped rank-sum tests, and with Savage and Sethuraman's [50] invariant SPRT of $H_0 : F = G$ versus a single Lehmann alternative $H_1 : F = G^\delta$ in the two-sample problem associated with population distribution functions $F$ and $G$. It was soon recognized that natural hypotheses in conventional (i.e., fixed sample size) nonparametrics could not be handled by invariant SPRTs, e.g., $F = G^\delta$ with fixed $\delta \neq 1$ is highly artificial. An alternative approach is to make use of the weak convergence of the sequence of suitably normalized nonparametric test statistics under the null hypothesis and under contiguous alternatives to Brownian motion, with drift 0 under $H_0$ and drift $\theta$ under a contiguous alternative, so that the nonparametric testing problem is asymptotically equivalent to that of testing whether the drift of Brownian motion is 0 or $\theta$. Sen and Ghosh [57], Sen [53, 54], Ghosh and Sen [18], Hall and Loynes [24] and Lai [34, 35] were some of the earliest applications of this approach. Chapters 9 and 11 of Sen's [56] monograph give an overview of these and subsequent developments in nonparametric sequential testing. Chapter 10 reviews related developments in nonparametric sequential interval and point estimation, while Chapters 2–8 provide the basic weak convergence theory.

Besides marking the publication of Sen's monograph, the year 1981 also witnessed another important event in sequential nonparametrics and its applications to clinical trials. The early termination of the Beta-Blocker Heart Attack Trial (BHAT) in October 1981 prior to its prescheduled end nine months later drew immediate attention of the medical community to the benefits of sequential methods

---

*Supported by the National Science Foundation and the National Institutes of Health.

[1]Stanford University, California, USA, e-mail: lait@stat.stanford.edu

[2]Genentech Inc, California, USA, e-mail: su.zheng@gene.com

and spurred important advances and new developments in sequential nonparametrics and semiparametrics in the following two decades, which are reviewed in Section 2. To tackle the complexity of time-sequential clinical trials with staggered patient entry and failure-time endpoints, new weak convergence results and asymptotic expansions were developed by making use of powerful tools involving continuous-time martingales, stochastic integrals and empirical process theory. Section 3 summarizes these and other results under the rubric of functional central limit theorems (CLT), also called "invariance principles" as in Sen's monograph, in sequential nonparametrics and semiparametrics.

In applying the asymptotic theory to the design and analysis of time-sequential trials with failure-time endpoints, several implementation issues need to be addressed and are discussed in Section 4. For example, the asymptotic approximations may be inadequate for the sample size considered, or may be difficult to compute directly. Another long-standing problem is related to terminal analysis of the trial. How can valid nonparametric/semiparametric confidence intervals be constructed for the primary and secondary endpoints when early stopping may occur during interim analysis of the trial? This problem is addressed in Section 4 which makes use of recent developments in Monte Carlo and resampling methods to resolve difficulties in the design and analysis of complex time-sequential trials.

## 2. Time-sequential rank tests in clinical trials

In a typical clinical trial to compare times to failure between two treatment groups $X$ and $Y$, $n$ patients enter the trial at different times during an accrual period, are randomly assigned to treatment $X$ or $Y$ and are then followed until they fail or withdraw from the study or until the study is terminated. In particular, BHAT was a multicenter, double-blind, randomized, placebo-controlled clinical trial designed to test the efficacy of long-term therapy with propranolol given to survivors of an acute myocardial infarction. The trial was scheduled for 4 years, with reviews of the data by a Data and Safety Monitoring Board (DSMB) at 11, 16, 21, 28, 34, 40 and 48 months. The trial design assumed an accrual rate of 149 patients per month for a period of 27 months, so the planned total number of patients was 4123. Another assumption is that each patient is randomized to placebo or treatment upon entering the trial, and is followed for a maximum of 3 years. The actual recruitment period was 27 months, within which 3837 patients were accrued from 136 coronary care units in 31 clinical centers, with 1916 patients randomized into the propranolol group and 1921 into the placebo group. Successive values of the standardized logrank statistics (see the next paragraph) at 11, 16, 21, 28, 34 and 40 months when the DSMB met were 1.68, 2.24, 2.37, 2.30, 2.34 and 2.82, respectively. Instead of continuing the trial to its scheduled end at 48 months, the DSMB recommended terminating it in their last meeting because of conclusive evidence in favor of propranolol. To adjust for repeated significance testing, the DSMB used critical values provided by O'Brien and Fleming's [43] method for normal observations; see Section 2.2. Since logrank statistics (rather than normal means) were actually used, the DSMB also appealed to Tsiatis' [62] result on the joint asymptotic normality of time-sequential logrank statistics.

To describe a typical time-sequential trial with staggered patient entry, we first introduce the following notation. Let $T_i' \geq 0$ denote the entry time and $X_i > 0$ the survival time (or time to failure) after entry of the $i$th subject in treatment group $X$, and let $T_j''$ and $Y_j$ denote the entry time and survival time after entry of the $j$th subject in treatment group $Y$, $1 \le i \le n'$, $1 \le j \le n''$. Let $n = n' + n''$. Thus the data at calendar time $t$ consist of $(X_i(t), \delta_i'(t)), i = 1, \ldots, n'$, and $(Y_j(t), \delta_j''(t)), j = 1, \ldots, n''$, where $X_i(t) = X_i \wedge \xi_i' \wedge (t - T_i')^+, Y_j(t) = Y_j \wedge \xi_j'' \wedge (t - T_j'')^+, \delta_i'(t) = I_{\{X_i(t) = X_i\}}, \delta_j''(t) = I_{\{Y_j(t) = Y_j\}}, m_{n,t}'(s) = \sum_{i=1}^{n'} I_{\{X_i(t) \ge s\}}, m_{n,t}''(s) = \sum_{j=1}^{n''} I_{\{Y_j(t) \ge s\}}$, and $\xi_i'(\xi_j'')$ denotes the withdrawal time, possibly infinite, of the $i$th ($j$th) subject in treatment group $X(Y)$. At a given calendar time $t$, one can compute, on the basis of the observed data from the two treatment groups, a rank statistic of the general form considered by Tsiatis [63]:

$$S_n(t) = \sum_{i=1}^{n'} \delta_i'(t) Q_n(t, X_i(t)) \left\{ 1 - \frac{m_{n,t}'(X_i(t))}{m_{n,t}'(X_i(t)) + m_{n,t}''(X_i(t))} \right\}$$

$$- \sum_{j=1}^{n''} \delta_j''(t) Q_n(t, Y_j(t)) \frac{m_{n,t}'(Y_j(t))}{m_{n,t}'(Y_j(t)) + m_{n,t}''(Y_j(t))}, \tag{1}$$

where $Q_n(t, s)$ is some weight function satisfying certain measurability assumptions. Let $H_{n,t}$ denote a product-limit-type estimator of the common distribution function of the two treatment groups under the null hypothesis, based on $\{(X_i(t), \delta_i(t), Y_j(t), \delta_j(t)) : i \le n', j \le n''\}$. Note that this setting is considerably more complicated than the progressively censored case considered by Chatterjee and Sen [9] and Sen [55, 56] in the context of life testing experiments. As will be reviewed in Section 3, for a general weight function of the form $Q_n(t, s) = \psi(H_{n,t}(s))$ in (1), $\{S_n(t)/\sqrt{n}, t \ge 0\}$ converges weakly to a Gaussian process with independent increments and variance function $V(t)$ under the null hypothesis, and contiguous alternatives. Gu, Lai and Lan [22] have pointed out that these time-sequential rank statistics based on censored data are natural extensions of classical rank statistics in fixed sample size (FSS) tests of $H_0 : F = G$, where $F$ is the common distribution function of the $X_i$ and $G$ is that of the $Y_j$, in which $X_i$ and $Y_j$ are completely observable (i.e., there is no censoring). Letting $R_i$ denote the rank of $X_i$ in the combined sample $\{X_1, \ldots, X_{n'}, Y_1, \ldots, Y_{n''}\}$, a classical rank statistic has the form $\sum_{i=1}^{n'} \varphi(R_i/n)$, where $\varphi : (0, 1] \to \mathcal{R}$ is a score function, and its extension to censored time-sequential data has the form (1) with $Q_n(t, s) = \psi(H_{n,t}(s))$, with $\psi$ related to $\varphi$ via the relation

$$\psi(u) = \varphi(u) - (1 - u)^{-1} \int_u^1 \varphi(t) \, dt, \qquad 0 < u < 1; \tag{2}$$

see Gu, Lai and Lan [22]. Taking $\psi(u) = (1 - u)^\rho$ ($\rho \ge 0$) yields the $G^\rho$ statistics proposed by Harrington and Fleming [25]. The case $\rho = 0$ corresponds to Mantel's [42] logrank statistic whose corresponding $\varphi$ is the asymptotically optimal score function for testing Lehmann alternatives. The case $\rho = 1$ corresponds to the generalization of Wilcoxon's statistic by Peto and Peto [44] and Prentice [48]. In the remainder of Section 2 it will be assumed that $Q_n(t, s) = \psi(H_{n,t}(s))$, unless stated otherwise.

The mean function of the limiting Gaussian process associated with $S_n(t)/\sqrt{n}$ is 0 under the null hypothesis $H_0 : F = G$ and is of the form $\mu_g(t)$ under contiguous alternatives that satisfy

$$\int_0^{t^*} \left| \frac{d\Lambda_G}{d\Lambda_F} - 1 \right| d\Lambda_F = O\left(\frac{1}{\sqrt{n}}\right), \quad \sqrt{n}\left\{\frac{d\Lambda_G}{d\Lambda_F}(s) - 1\right\} \to g(s) \tag{3}$$

as $n \to \infty$, uniformly over closed subintervals of $\{s \in [0, t^*] : F(s) < 1\}$, where $\Lambda_F$ and $\Lambda_G$ are the cumulative hazard functions of $F$ and $G$; see Section 3. Moreover, $\mu_g(t) = V(t)$ when $\psi$ is an asymptotically optimal score function of the form $\psi(\cdot) = g(F^{-1}(\cdot))$. In practice, the actual alternatives are unknown and $\mu_g$ need not even be monotone when $\psi$ is not optimal for the actual alternatives, such as using logrank statistics for non-proportional hazards (i.e., not Lehmann) alternatives, which allow repeated significance tests based on $S_n(t)$ to achieve *both* savings in study duration and increase in power over the fixed-duration test based on $S_n(t^*)$; see Example 1 of Gu and Lai [20].

## 2.1. Estimation of $V(t)$ and two time scales

An estimate of the variance of $S_n(t)$ under $H_0$ can be used to estimate the variance $V(t)$ of the limiting Gaussian approximation to $S_n(t)/\sqrt{n}$ under $H_0$ and under contiguous alternatives. Two commonly used variance estimates are

$$V_n(t) = \int_0^t \frac{\psi^2(H_{n,t}(s)) m'_{n,t}(s) m''_{n,t}(s)}{(m'_{n,t}(s) + m''_{n,t}(s))^2} \, d(N'_{n,t}(s) + N''_{n,t}(s)), \qquad (4a)$$

$$V_n(t) = \int_0^t \frac{\psi^2(H_{n,t}(s))}{(m'_{n,t}(s) + m''_{n,t}(s))^2} \{(m''_{n,t}(s))^2 \, dN'_{n,t}(s) + (m'_{n,t}(s))^2 \, dN''_{n,t}(s)\}, \quad (4b)$$

where $N'_{n,t} = \sum_{i=1}^{n'} I_{\{X_i \le \xi'_i \wedge (t-T'_i)^+ \wedge s\}}$, $N''_{n,t} = \sum_{j=1}^{n''} I_{\{Y_j \le \xi''_j \wedge (t-T''_j)^+ \wedge s\}}$. As a compromise between these two choices, Gu and Lai [19] also considered

$$V_n(t) = \{(4a) + (4b)\}/2. \qquad (4c)$$

For all three estimates, $n^{-1} V_n(t)$ converges in probability to $V(t)$ under $H_0$ and under contiguous alternatives. Hence, letting $v = n^{-1} V_n(t)$ and $W(v) = n^{-1/2} S_n(t)$, we can regard $W(v), v \ge 0$, as the standard Wiener process under $H_0$. Moreover, if $\psi$ is a scalar multiple of the asymptotically optimal score function, then we can also regard $W(v), v \ge 0$, as a Wiener process with some drift coefficient under contiguous alternatives.

As pointed out by Lan and DeMets [39], there are two time scales in time-sequential trials with failure-time endpoints. One is *calendar time* $t$ and the other is *information time* $V_n(t)$, and there is no simple relation between them. The calendar times of interest are those when interim (including final) analyses are performed and they are usually specified, at least approximately, in the trial protocol. The information time $V_n(t)$ is typically unknown before time $t$ unless restrictive assumptions are made *a priori*. The two time scales have led to several long-standing difficulties in applying the asymptotic theory to repeated significance tests based on time-sequential rank statistics and to interval estimation following these tests. These difficulties have been resolved in the past five years, as reviewed in Sections 2.2 and 4.2.

## 2.2. Stopping boundaries for time-sequential rank tests

Let $0 < t_1 < \cdots < t_k = t^*$ be prescribed times for periodic reviews of the data. To test $H_0 : F = G$ with the time-sequential rank statistics $S_n(t_i)$, Slud and Wei [61] introduced the following simple approach. First choose positive numbers $\alpha_1, \ldots, \alpha_k$ such that $\sum_1^k \alpha_i = \alpha$ (= the overall significance level). Then use the multivariate normal approximation to the null distribution of $(S_n(t_i)/V_n^{1/2}(t_i))_{1 \le i \le k}$ to determine $d_1, \ldots, d_k$ recursively by

$$P_{H_0}\{|S_n(t_j)|/V_n^{1/2}(t_j) \ge d_j \text{ and } |S_n(t_i)|/V_n^{1/2}(t_i) < d_i \text{ for all } i < j\} = \alpha_j. \quad (5)$$

With the $d_j$ thus determined, the Slud-Wei repeated significance test rejects $H_0$ whenever $|S_n(t_j)| \ge d_j V_n^{1/2}(t_j)(1 \le j \le k)$ and stops the trial at the first $t_j$ this occurs (or at $t^*$ if this does not occur for $1 \le j \le k$).

The Slud-Wei method does not provide practical guidelines concerning how the $\alpha_j$ in (5) should be chosen. Lan and DeMets [38] and Lan, DeMets and Halperin [40] proposed to derive the $\alpha_j$ from an *error spending function*, which specifies how fast we can spend the Type I error $\alpha$ over time. To begin with, let $\{W(v), 0 \le v \le 1\}$ be the standard Wiener process and consider the stopping rule $T = \inf\{v \in [0,1] : |W(v)| \ge h(v)\}(\inf \emptyset = \infty)$, where $h$ is a positive function on $[0,1]$ such that $P\{T = 0\} = 0$ and $P\{T \le 1\} = \alpha$. The error spending function is $A(v) = P\{T \le v\}, 0 \le v \le 1$. Taking $v$ to represent the proportion of information accumulated at time $t$, with $A(0) = 0$ and $A(1) = \alpha$. In particular, suppose that instead of survival data one has immediate responses from the patients who enter the study serially and are randomized to either treatment, with a target sample size of $n$ at the scheduled end of the trial. Lan and DeMets [39] call such trials "maximum information trials". Here the proportion of information accumulated at time $t_i$ of interim analysis $v_i = n_i/n$, where $n_i$ is the total number of patients available at $t_i$. Hence Lan and DeMets [38] proposed to choose $\alpha_j = A(v_j) - A(v_{j-1})$ in (5). For the time-sequential rank statistics (1) in what Lan and DeMets [39] call "maximum duration trials", the asymptotic null variance of $S_n(t_i)$ is no longer proportional to the sample size $n_i$ at $t_i$ and a natural analogue of $n_i/n$ here is $V_n(t_i)/V_n(t^*)$.

Let $Z_1, Z_2, \ldots$ be i.i.d. normal random variables with unknown mean $\theta$ and known variance 1. To test $H : \theta = 0$ at level $\alpha$, the Neyman-Pearson test rejects $H$ if $|\sum_{i=1}^{k} Z_i| \ge z_\alpha \sqrt{k}$, where $1 - \Phi(z_\alpha) = \alpha$. Sample size calculation in clinical trial applications typically assume an alternative $\theta$ of particular interest and find the $k$ that attains some given power $1 - \beta$ at $\theta$. The basic idea behind Haybittle's [26] repeated significance test is to keep $k$ and $\alpha$ as the maximum sample size and significance level but to allow for early stopping when the data are monitored sequentially, at the expense of some minor loss in power at $\theta$. This leads to the stopping rule $\tau_b = \min(k, \inf\{n \ge 1 : |\sum_{i=1}^{n} Z_i| \ge b\sqrt{n}\})$ and terminal decision rule that rejects $H$ if $\tau_b < k$ or if $\tau_b = k$ and $|\sum_{i=1}^{k} Z_i| \ge c\sqrt{k}$. Since we require the loss in power at $\theta$ to be small relative to the fixed sample size test and also require the maximum sample size to be the same as the fixed sample size $k$, it is clear that $c$ has to be near $z_\alpha$, implying that $P_0(\tau_b < k)$ is small in comparison with $\alpha$. In particular, the Haybittle-Peto method [26, 45] in the field of clinical trials uses some relatively large value of $b$, such as 3, and conventional critical values of $c$ for the final test when the number $k$ of interim analyses is small. The fact that $P_0(\tau_b < k)$ is typically small relative to $\alpha$ (or equivalently that most of the Type I error is to be spent at the terminal data $t^*$) suggests that using an elaborate Lan-DeMets boundary determination procedure would not lead to substantial improvement over a simple procedure of the Haybittle-Peto type. Lai and Shih [37] recently developed a theory of group sequential tests, based on $Z_1, Z_2, \ldots$, for the parameter $\theta$ of an exponential family of densities $f_\theta(z) = e^{\theta z - \psi(\theta)}$ with respect to some measure on the real line. This theory yields the following modified Haybittle-Peto test as an asymptotically optimal solution. Let $S_n = Z_1 + \cdots + Z_n$. To test the one-sided hypothesis $H_0 : \theta \le \theta_0$ at significance level $\alpha$, suppose no more than $M$ observations are to be taken. The fixed-sample-size test that rejects $H_0$ if $S_M \geq c_\alpha$ has maximal power at any alternative $\theta > \theta_0$, in particular at the alternative $\theta(M)$ "implied" by $M$ (in the sense that $M$ can be derived from the assumption that the above fixed-sample-size test has some prescribed power $1 - \tilde{\alpha}$ at $\theta(M)$). Although the protocol of a clinical trial typically justifies its choice of sample size by stating some conventional level (such as 80% or 90%) at a specified alternative, one often does not have much information prior to the trial to come up with a realistic alternative. Under the constraint of $M$ on the sample size, it is desirable to adapt to the information on the actual $\theta$ gathered during the course of the trial, allowing early stopping at times of interim analysis so that the test has nearly optimal power and expected sample size properties.

To achieve these goals in a group sequential test with $k$ groups and group sizes $n_1, n_2 - n_1, \ldots, n_k - n_{k-1}$ so that $n_k = M$, Lai and Shih [37] use a rejection region of the form $S_{n_k} \geq c$ at the $k^{th}$ analysis, where $c > c_\alpha$ but $c$ does not differ much from $c_\alpha$. For the first $k - 1$ analyses, they use a stopping region of the form

$$\widehat{\theta}_{n_i} > \theta_0 \quad \text{and} \quad n_i I(\widehat{\theta}_{n_i}, \theta_0) \geq b, \text{ or} \tag{6a}$$

$$\widehat{\theta}_{n_i} < \theta(M) \quad \text{and} \quad n_i I(\widehat{\theta}_{n_i}, \theta(M)) \geq \widetilde{b}, \tag{6b}$$

for $1 \leq i \leq k - 1$, where $I(\theta, \lambda)$ denotes the Kullback-Leibler information number $E_\theta\{\log[f_\theta(Z_i)/f_\lambda(Z_i)]\}$. If (6a) holds, reject $H_0$ upon stopping. If stopping occurs with (6b), accept $H_0$. In case stopping does not occur in the first $k - 1$ analyses, reject $H_0$ if $S_{n_k} \geq c$. The thresholds $b, \widetilde{b}$ and $c$ are so chosen that $P_{\theta_0}$ (Test rejects $H_0$)$= \alpha$ and the power of the test at $\theta(M)$ does not differ much from its upper bound $1 - \tilde{\alpha}$. In the special case of $\theta_0 = 0$ and normal $Z_i$ with mean $\theta$ and variance 1, $\psi(\theta) = \theta^2/2$ and $I(\theta, \lambda) = (\theta - \lambda)^2/2$, so the test reduces to that in the preceding paragraph. Previous works on efficient group sequential tests have used the expected sample size at an alternative, or more generally a weighted average of expected sample sizes over a set of parameter values, as the optimization criterion while controlling the error probabilities under the null hypothesis and a specified alternative at prescribed levels; see Pocock [47], Wang and Tsiatis [67], Kim and DeMets [32], Eales and Jennison [15] and Barber and Jennison [4]. There are several practical difficulties with this approach to efficient group sequential design. First, even though the mean of the random sample size is minimized at some alternative, the maximum sample size can be substantially larger than the mean and also the fixed sample size. Secondly, the optimization problem requires precise specification of the relative sizes of all groups, e.g. equal group sizes, but it is often not feasible to do so prior to the trial because interim analyses are usually scheduled at calender times for administrative reasons. Thirdly, it may be difficult to come up with a realistic alternative before data are collected from the trial, but the optimization problem depends on the chosen alternative. Clearly efficiency of a group sequential test depends not only on the choice of the stopping rule but also on the test statistics used. To decouple these two issues, Lai and Shih [37] consider the one-parameter exponential family, for which sufficient statistics are the sample means that are maximum likelihood estimators of $\psi'(\theta)$. Although the normal case is usually chosen to be the prototype in the group sequential literature, Lai and Shih choose the exponential family because linearity of $\psi'$ in the normal case obscures the general form of (nearly) optimal test statistics and stopping boundaries. Making use of Hoeffding's [27] lower bound for the expected sample size of a test that has Type I error probability $\alpha$ at $\theta_0$ and Type II error probability $\tilde{\alpha}$ at $\theta(M)$, they establish the asymptotic efficiency of the modified Haybittle-Peto test by showing that the expected sample size of the test under $P_\theta$ attains Hoeffding's lower bound asymptotically as $\alpha \to 0$ and $\bar{\alpha} \to 0$. They have also carried out extensive simulation studies of the expected sample size and power function of the modified Haybittle-Peto test in the $N(\theta, 1)$ case, comparing it with different classes of group sequential tests that have been developed primarily for this normal setting in the literature. Their numerical results demonstrate the advantages of the modified Haybittle-Peto test, which is flexible and efficient and can "self-tune" to the unknown parameters during the course of the trial, under prespecified constraints on the maximum sample size and Type I error probability.

The preceding modified Haybittle-Peto test has in fact been proposed earlier by Gu and Lai [20] for repeated significance testing, which involves a maximum of $k$ significance tests, based on the time-sequential censored rank statistics (1). Gu and Lai proposed to determine $b$ such that $P_0(\tau_b < k) = \epsilon\alpha$, where $0 < \epsilon < 1$ is small and $\tau_b$ is the stopping rule defined above for normal $Z_i$. With $b$ thus chosen, their repeated significance test stops the trial and rejects $H_0$ at $t_i < t^*$ if $|S_n(t_i)| \geq bV_n^{1/2}(t_i)$. If the trial proceeds to the terminal date $t^*$, their test rejects $H_0$ if $|S_n(t^*)| \geq cV_n^{1/2}(t^*)$, where $c$ is so chosen that

$$P\{|W(V_n(t_k))| \geq cV_n^{1/2}(t_k) \text{ or } |W(V_n(t_i))| \geq bV_n^{1/2}(t_i)$$

$$\text{for some } i < k | V_n(t_1), \ldots, V_n(t_k)\} = \alpha, \qquad (7)$$

in which $t_k = t^*$ and $\{W(v), v \geq 0\}$ is a standard Wiener process independent of $\{(X_i, \xi_i', t_i', Y_i, \xi_i'', T_i''), i \geq 1\}$. Letting $a_j = V_n(t_j)$ and $d_k = c, d_j = b$ for $1 \leq j \leq k - 1$, the probability (7) can be written as a sum of the probabilities $P\{|W(a_j)| \geq d_j\sqrt{a_j}$ and $W(a_i) < d_i\sqrt{a_i}$ for all $i < j\}$, which can be computed by the recursive numerical integration algorithm of Armitage, McPherson and Rowe [3]. The choice of $c$ in (7) is not predetermined at the beginning of the trial but depends on the actual values of $V_n(t_1), \ldots, V_n(t_k)$, allowing great flexibility in how information accumulates at different times of interim analyses. This modified Haybittle-Peto time-sequential test based on (1), with $t$ restricted to the set $\{t_1, \ldots, t_k\}$ of calendar times at which interim analyses are conducted, yields an efficient stopping rule that circumvents the difficulty of "calendar time" versus "information time" in the error spending approach, which Scharfstein and Tsiatis [51] proposed to address by using simulations at each interim analysis to estimate the maximum information under the null hypothesis.

## 2.3. Adjustments for other covariates in testing treatment effects

It is widely recognized that tests of treatment effects based on the rank statistics (1) may lose substantial power when the effects of other covariates are strong. In nonsequential trials, a commonly used method to remedy this when logrank statistics are used is to assume the proportional hazards regression model and to use Cox's partial likelihood approach to adjust for other covariates. Tsiatis, Rosner and Tritchler [65], Gu and Ying [23] and Bilias, Gu and Ying [5] have developed group sequential tests using this approach. Instead of relying on the proportional hazards model to adjust for concomitant variables, it is useful to have other methods for covariate adjustment, especially in situations where other score functions than the logrank are used in (1) to allow for the possibility of non-proportional hazards alternatives. Lin [41] and Gu and Lai [20] have proposed alternative covariate adjustment methods based on rank estimators and $M$-estimators in accelerated failure time models. In the absence of more efficient computational schemes, these estimators took much longer to compute than those based on the proportional hazards model. The situation has subsequently improved with the algorithm of Jin et al. [29] for rank estimators and Kim and Lai [33] for M-estimators.

## 3. Functional CLT in sequential non-(or semi-)parametrics

Let $\widehat{F}_n$ be the empirical distribution function of i.i.d. random variables $X_1, X_2, \ldots$ with common continuous distribution function $F$, and let $F_U(u) = u$ for $0 \le u \le 1$, the distribution function of a uniform random variable. Whereas the weak convergence of $\sqrt{n}(\widehat{F}_n \circ F^{-1} - F_U)$ to Brownian bridge, when applied in conjunction with the functional delta method, provides a basic tool for deriving asymptotic normality of nonparametric/semiparametric statistics, sequential nonparametrics and semiparametrics involve an additional time parameter, which results in weak convergence to Gaussian random fields (i.e., multiparameter processes). Sen's [56] monograph therefore focuses on weak convergence in $C([0, 1]^k)$ or $D([0, 1]^k)$ with $k \ge 2$. In particular, $[nt](\widehat{F}_{[nt]} \circ F^{-1} - F_U)/\sqrt{n}$ converges weakly in $D([0, 1]^2)$ to a Kiefer process $K(t, s)$ whose covariance function is given by

$$\text{Cov}(K(t, s), K(t', s')) = (t \wedge t')(s \wedge s' - ss'),$$

i.e., $K(\cdot, s)$ is Brownian bridge for every fixed $s$ and $K(t, \cdot)$ is Brownian motion for every fixed $t$. An alternative approach is to represent the nonparametric statistics (e.g. $U$-statistics, rank statistics and linear combinations of order statistics) in terms of partial sums $S_n$ of i.i.d. random variables plus negligible remainders; see Lai [34, 35].

When the $X_i$ are not fully observable because of censoring by $\xi_i$, it is more convenient to work with the cumulative hazard function $\Lambda(x) = \int_{-\infty}^{x}(1 - F(s-))^{-1} dF(s)$, which can be consistently estimated by the nonparametric maximum likelihood estimator $\widehat{\Lambda}_n(x) = \int_{-\infty}^{x}(m_n(s))^{-1} dN_n(s)$, where

$$m_n(s) = \sum_{i=1}^{n} I_{\{X_i \wedge \xi_i \ge s\}}, \quad N_n(s) = \sum_{i=1}^{n} I_{\{X_i \le \xi_i \wedge s\}}.$$

Making use of the key property that $\{N_n(s) - \int_{-\infty}^{s} m_n(t) \, d\Lambda(t), s \ge 0\}$ is a square integrable martingale with predictable variation process $\int_{-\infty}^{s} m_n(t) \, d\Lambda(t)$, we can apply martingale central limit theorems to prove asymptotic normality of statistics of the form $\int_{-\infty}^{x} Q_n(s) \, d\widehat{\Lambda}_n(s)$, where $Q_n(\cdot)$ is a predictable process. The monograph by Andersen et al. [1] summarizes this martingale approach to functional central limit theorems for nonparametric and semiparametric statistics based on censored data and their applications.

This martingale approach can be easily extended to the vector of the time-sequential rank statistics $(S_n(t_1), \ldots, S_n(t_k))/\sqrt{n}$ involving a fixed set of calendar times $t_1, \ldots, t_k$, where $S_n(t)$ is defined in (1), as shown by Tsiatis [62, 63]. However, proving weak convergence of the continuous-time process $\{n^{-1/2}S_n(t), 0 \le t \le t^*\}$, where $t^*$ is the terminal tie of the study, is much more difficult even though it involves only checking an additional tightness condition. By combining certain maximal inequalities for continuous-time martingales with empirical process theory, Gu and Lai ([19], Lemma 2 and Appendix) established the desired tightness under some weak conditions on the weights $Q_n(t, s)$ in (1) and the assumptions that

$$b'(t, s) = \lim_{m \to \infty} m^{-1} \sum_{i=1}^m P\{\xi_i' \geq s, t - T_i' \geq s\},$$

$$b''(t, s) = \lim_{m \to \infty} m^{-1} \sum_{j=1}^m P\{\xi_j'' \geq s, t - T_j'' \geq s\}$$

exist and are continuous in $0 \leq s \leq t$, that

$$n'/n \to \gamma \text{ as } n \to \infty \text{ with } 0 < \gamma < 1,$$

and that $F$ and $G$ are continuous. Instead of $S_n(t)$, they considered $S_n(t, s)$ which is defined by (1) but with $\sum_{1 \leq i \leq n'}$ and $\sum_{1 \leq j \leq n''}$ replaced by $\sum_{i:X_i(t) \leq s}$ and $\sum_{j:Y_j(t) \leq s}$. In particular, for the case $Q_n(t, s) = \psi(H_n(t, s))$ such that $(1-x)^\beta \psi(x)$ is a function of bounded variation on $[0, 1]$ for some $0 \leq \beta < \frac{1}{2}$, their weak convergence theory for the random field $\{S_n(t, s), 0 \leq t \leq t^*, 0 \leq s \leq t^*\}$ yields the following results for $S_n(t)$:

(i) For fixed $F$ and $G$, $\{n^{-1/2}(S_n(t) - \mu_n(t)), 0 \leq t \leq t^*\}$ converges weakly in $D[0, t^*]$ to a zero-mean Gaussian process and $n^{-1}\mu_n(t)$ converges in probability as $n \to \infty$, where

$$\mu_n(t) = \int_0^t \psi(H_{n,t}(s)) \frac{m'_{n,t}(s) m''_{n,t}(s)}{m'_{n,t}(s) + m''_{n,t}(s)} (d\Lambda_F(s) - d\Lambda_G(s)).$$

(ii) Let $\{Z(t), 0 \leq t \leq t^*\}$ denote the zero-mean Gaussian process in (i) when $F = G$. This Gaussian process has independent increments and

$$\text{Var}(Z(t)) = \gamma(1-\gamma) \int_0^t \frac{\psi^2(F(s)) b'(t, s) b''(t, s)}{\gamma b'(t, s) + (1-\gamma) b''(t, s)} dF(s).$$

(iii) For fixed $F$ (and therefore $\Lambda_F$ also), suppose that as $n \to \infty$, $G \to F$ such that (3) holds uniformly over closed subintervals $I$ of $\{s \in [0, t^*] : F(s) < 1\}$ and $\sup_{s \in I} |g(s)| < \infty$. Then $\{n^{-1/2} S_n(t), 0 \leq t \leq t^*\}$ converges weakly in $D[0, t^*]$ to $\{Z(t) + \mu(t), 0 \leq t \leq t^*\}$, where $Z(t)$ is the same Gaussian process as that in (ii) and

$$\mu(t) = -\gamma(1-\gamma) \int_0^t \frac{\psi(F(u)) g(u) b'(t, u) b''(t, u)}{\gamma b'(t, u) + (1-\gamma) b''(t, u)} dF(u).$$

From (ii) and (iii), the limiting Gaussian process of $\{n^{-1/2} S_n(t), t \geq 0\}$ has independent increments under $H_0 : F = G$ and under contiguous alternatives.

Previous results in the literature only treated the case $F = G$. In particular, assuming the $T_i'$ (and $T_i''$, $\xi_i'$, $\xi_i''$, respectively) to be i.i.d., Tsiatis [63] showed that $(n^{-1/2} S_n(t_1), \ldots, n^{-1/2} S_n(t_k))$ has a limiting multivariate normal distribution for any $k$. Sellke and Siegmund [52] proved tightness and weak convergence in the case where the $S_n(t)$ are logrank statistics, and more generally where $S_n(t)$ is the score statistic in the proportional hazards model, without assuming the $T_i'$ (or $T_i''$, $\xi_i'$, $\xi_i''$) to be i.i.d.. Slud [60] considered weighted logrank statistics with weights that do not depend on $t$.

## 4. Implementation issues in sequential non-(or semi-)parametrics

This section considers certain implementation issues in sequential nonparametrics and semiparametrics, particularly in the context of design and analysis of complex clinical trials with failure-time endpoints and periodic data reviews. Although the asymptotic joint normality of the sequential non-(or semi-)parametric statistics greatly simplifies their seemingly intractable distributions, the adequacy of these approximations may be questionable. For example, using the normal approximation to evaluate the power of a sequential rank test at an alternative that is assumed to be "contiguous" may be unreliable as it is difficult to assess whether the alternative is actually contiguous for the sample size associated with a prescribed stopping rule that allows early termination. Another important issue is related to the construction of confidence intervals following sequential tests. The normal approximation can be applied if the stopped information time is asymptotically nonrandom by Anscombe's [2] theorem, but it is difficult to assess whether that is indeed the case in practice. The two time scales in time-sequential clinical trials further increase the difficulty. In this section we use recent developments in Monte Carlo and resampling methods to address these issues.

### 4.1. Direct Monte Carlo and its enhancements

The Monte Carlo simulation method provides a flexible and practical way to compute the power and expected duration of time-sequential tests, and also to check the adequacy of the normal approximation to the type I error probability under various scenarios of baseline survival, censoring pattern, noncompliance, and accrual rate. To provide the clinical trial designer with a tool to perform these Monte Carlo simulations, Gu and Lai [21] developed a simulation program which gives the user some options for choosing the stopping boundary, including the modified Haybittle-Peto-type boundary. They also incorporated this power calculation program into another program that computes the sample size of a group sequential trial having a prescribed power at given baseline and alternative distributions. Adjustment for other covariates in time-sequential tests of treatment effects, however, is unavailable in the program developed by Gu and Lai [21]. Because of the computational complexity of time-sequential tests and because the Monte Carlo simulations used to compute power and type I error probability should not take too long to run for the software to be "user-friendly", the direct Monte Carlo approach used by Gu and Lai [21], which is already slow to run a "bare-bones" trial design, cannot absorb the additional computational costs of covariate adjustment without further slowing it down substantially.

Importance sampling and related exponential tilting techniques are powerful methods for Monte Carlo evaluation of small probabilities of events observed up to a stopping time. The basic idea comes from the likelihood ratio identity $P(F \cap \{T < \infty\}) = E_Q(L_T 1_{F \cap \{T < \infty\}})$ for all $F \in \mathcal{F}_T$, where $T$ is a stopping time and $L_n$ is the likelihood ratio; see Siegmund [59], page 13. For complicated problems, however, implementation of these importance sampling methods may be difficult due to difficulties of sampling from $Q$. Instead of sampling directly from $Q$, we can generate $B$ sequences from a more convenient distribution $\widetilde{Q}$ and then use a resampling step to convert samples from $\widetilde{Q}$ to samples from $Q$. There are results from the bootstrap literature that shed light on how to carry this out; see Johns [30], Davison and Hinkley [13], Do and Hall [14] and Hu and Su [28]. The following example illustrates this idea with the choice of resampling weights for bootstrap estimation of the sampling distribution of a two-sample Mann-Whitney statistic from unknown $(F, G)$.

**Example 1.** Denote the two independent samples by $\mathcal{X} = \{X_1, \ldots, X_m\}$ and $\mathcal{Y} = \{Y_1, \ldots, Y_n\}$ drawn from $F$ and $G$, respectively. The Mann-Whitney statistic is given by $U = \sum_{i=1}^{m} \sum_{j=1}^{n} U(i, j)$, where

$$
U(i, j) = \begin{cases} +1 & \text{if } X_i > Y_j, \\ 0 & \text{if } X_i = Y_j, \\ -1 & \text{if } X_i < Y_j. \end{cases}
$$

To estimate the probability $P(U \leq x)$ from $\mathcal{X}$ and $\mathcal{Y}$, the bootstrap method estimates $P(U^* \leq x | \mathcal{X}, \mathcal{Y})$ by Monte Carlo simulations, where $U^*$ is the value of the statistic calculated from the bootstrap resamples $\mathcal{X}^*$ and $\mathcal{Y}^*$. To reduce the Monte Carlo variance of the bootstrap approach, we fix the first sample and resample from the second sample with resampling weights $p_i(i = 1, \ldots, n)$ instead of the usual weights $1/n$ corresponding to the empirical distribution of $\mathcal{Y}$. We can estimate $P(U^* \leq x | \mathcal{X}, \mathcal{Y})$ by $\hat{\pi} := B^{-1} \sum_{b=1}^{B} U_b \prod_{i=1}^{n} (np_i)^{-M_{bi}}$, where $B$ is the number of resamples, $U_b$ is the value of the test statistic calculated from the $b$th resample, and $M_{bi}$ denotes the number of times that $Y_i$ appears in the $b$th resample. The optimal resampling weights are those that minimize $\text{Var}(\hat{\pi})$, or equivalently, minimize

$$
\begin{aligned}
& E\{I(U^* \leq x) \prod_{i=1}^{n} (np_i)^{-M_i^*} \mid \mathcal{X}, \mathcal{Y}\} \\
= \ & E\{I(\sum_{i=1}^{n} M_i^* u_i \leq x) \prod_{i=1}^{n} (np_i)^{-M_i^*} \mid \mathcal{X}, \mathcal{Y}\} \\
= \ & E\{I(\sum_{i=1}^{n} M_i^* (u_i - \bar{u}) \leq x - n\bar{u}) \prod_{i=1}^{n} (np_i)^{-M_i^*} \mid \mathcal{X}, \mathcal{Y}\} \\
= \ & E\{I(\sum_{i=1}^{n} M_i^* \frac{u_i - \bar{u}}{\sqrt{\sum (u_i - \bar{u})^2}} \leq \frac{x - n\bar{u}}{\sqrt{\sum (u_i - \bar{u})^2}}) \prod_{i=1}^{n} (np_i)^{-M_i^*} \mid \mathcal{X}, \mathcal{Y}\} \\
= \ & E\{I(\sum_{i=1}^{n} M_i^* \tilde{u}_i \leq \tilde{x}) \prod_{i=1}^{n} (np_i)^{-M_i^*} \mid \mathcal{X}, \mathcal{Y}\} \\
\sim \ & E\{I(N_1 \leq \tilde{x}) e^{N_2} \mid \mathcal{X}, \mathcal{Y}\},
\end{aligned}
\tag{8}
$$

where $(N_1, N_2)$ is bivariate normal random vector with mean $(0, \frac{1}{2}s^2)$, variances $1, s^2$ and covariance $\sum \tilde{u}_i \delta_i$. Here $M_i^*$ denotes the number of times $Y_i$ appears in a bootstrap resample and

$$
u_j = \sum_{i=1}^{m} U(i, j), \quad \delta_i = -\log(np_i), \quad s^2 = \sum \delta_i^2,
$$

$$
\bar{u} = \frac{\sum_{i=1}^{n} u_i}{n}, \quad \tilde{x} = \frac{x - n\bar{u}}{\sqrt{\sum (u_i - \bar{u})^2}}, \quad \tilde{u}_i = \frac{u_i - \bar{u}}{\sqrt{\sum (u_i - \bar{u})^2}}.
$$

Since $E\{I(N_1 \leq \tilde{x}) e^{N_2} \mid \mathcal{X}, \mathcal{Y}\} = \Phi(\tilde{x} - s\rho) e^{s^2}$, where $\rho = \sum u_i \delta_i / \sqrt{\sum u_j^2} s$ and $\Phi$ denotes the cumulative distribution function of the standard normal distribution,

TABLE 1
*Bootstrap estimates of tail probabilities of 2-sample Wilcoxon statistic and their normal approximations*

| $\Phi(\tilde{x})$ | Direct resampling | Importance resampling |
|---|---|---|
| 0.005 | $0.0044 \pm 0.0046$ | $0.0045 \pm 0.0006$ |
| 0.01 | $0.0096 \pm 0.0078$ | $0.0110 \pm 0.0012$ |
| 0.025 | $0.0279 \pm 0.0132$ | $0.0279 \pm 0.0030$ |
| 0.05 | $0.0472 \pm 0.0139$ | $0.0464 \pm 0.0040$ |
| 0.1 | $0.1023 \pm 0.0214$ | $0.0994 \pm 0.0078$ |

the resampling weights that minimize (8) are

$$p_i = \frac{e^{-A\tilde{u}_i}}{\sum_{j=1}^n e^{-A\tilde{u}_j}}, \ 1 \le i \le n,$$

where $A = A(x) > 0$ is chosen to minimize $\Phi(\tilde{x} - A)e^{A^2}$.

Table 1 gives the bootstrap estimate of $P(U^* \le x | \mathcal{X}, \mathcal{Y})$ and the standard errors for the usual resampling weights $1/n$ and for the optimal resampling weights given above. Here $m = 30, n = 25$, $B = 500$ bootstrap samples are used and $F = G$ is the exponential distribution with median 3. Table 1 shows that this importance resampling approach yields considerably smaller standard errors than the direct resampling approach. Further details of importance sampling with resampling for Monte Carlo computation of the power and Type I error of group sequential or time-sequential nonparametric/semiparametric tests will be given elsewhere.

### 4.2. Confidence intervals following time-sequential trials

Although group sequential or time-sequential tests are attractive in clinical trials because they allow for early termination while preserving the overall significance level of the test and can adapt to information gathered during the course of the trial, the use of a stopping rule may introduce substantial difficulties in constructing valid confidence intervals, which has inhibited the applications of sequential methodology. Siegmund [58] introduced an exact method, based on ordering the sample space $(T, S_T)$ in the following way, to construct exact confidence intervals for the mean of a normal population with known variance following a repeated significance test. Suppose $Z_i$ has variance 1 and $T$ is a two-sided stopping rule of the form $T = \inf\{n \in J : S_n \ge b_n \text{ or } S_n \le a_n\}$, where $S_n = Z_1 + \cdots + Z_n$ and $J$ is a finite set of positive integers. Siegmund ordered the sample space of $(T, S_T)$ as follows: $(t, s) > (t', s')$ whenever (i) $t = t'$ and $s > s'$, or (ii) $t < t'$ and $s \ge b_t$, or (iii) $t > t'$ and $s' \le a_{t'}$. Let $\mu_c$ denote the value of $\mu$ for which $P_\mu\{(T, S_T) \ge (t, s)_{\text{obs}}\} = c$, where $(t, s)_{\text{obs}}$ denotes the observed value of $(T, S_T)$. Siegmund's confidence interval is $\mu_\alpha \le \mu \le \mu_{1-\alpha}$, which has coverage probability $1 - 2\alpha$. Tsiatis, Rosner and Mehta [64] extended Siegmund's method to the group sequential tests of Pocock [46] and O'Brien and Fleming [43]. Alternative orderings of the sample space were subsequently introduced by Chang and O'Brien [8], Rosner and Tsiatis [49], Chang [7] and Emerson and Fleming [17]. To remove the normal assumption in these exact methods, a natural way is to extend Efron's [16] bootstrap confidence intervals from the fixed sample size to the group sequential setting. Chuang and Lai [10, 11], however, have shown that bootstrap confidence intervals following group sequential tests have inaccurate coverage probabilities because the approximate pivots for a fixed sample size $n$ may not remain to be approximate pivots when $n$ is replaced by a random stopping time $T$. By integrating the main ideas behind the exact and bootstrap methods, they proposed a *hybrid resampling method* and considered in particular those based on Siegmund's and other orderings for deriving the exact methods. They showed that hybrid resampling still works well in situations where the bootstrap method fails.

To extend the hybrid resampling method to more general statistics $\Psi_T$ and stopping rules $T$, e.g. those arising in time-sequential nonparametrics and semi-parametrics, Lai and Li [36] recently introduced the following ordering scheme for the sample space of a stochastic process $X_u$(in which $u$ denotes either discrete or continuous time) that is observed up to a stopping time $T$. Let $\Psi_t, t \leq T$, be real-valued statistics based on $\{X_t, t \leq T\}$. A total ordering of the sample space of $X$ can be defined via $(\Psi_t, t \leq T)$ as follows:

$$X \geq x \quad \text{if and only if} \quad \Psi_{T \wedge t} \geq \psi_{T \wedge t}, \tag{9}$$

in which $(\psi_s, s \leq t)$ is defined from $x = (x_s, s \leq t)$ in the same way as $(\Psi_s, s \leq t)$ is defined from $X$. In particular, they applied this ordering scheme to construct confidence intervals for the regression parameter $\beta$ in Cox's [12] proportional hazards model

$$P\{y \leq Y_i \leq y + dy | Y_i \geq y, z_i\} = e^{\beta z_i} d\Lambda(y)$$

with baseline cumulative hazard function $\Lambda$. Differentiation of the log partial likelihood at $\beta = 0$ (the null hypothesis) and calendar time $t$ yields Cox's score statistic

$$S_n(t) = \sum_{i=1}^{n} \delta_i(t) \left\{ z_i - \left( \sum_{j \in R_i(t)} z_j \right) \middle/ |R_i(t)| \right\}, \tag{10}$$

where $R_i(t) = \{j : Y_j(t) \geq Y_i(t)\}$ and $|R_i(t)|$ denotes the size of the "risk set" $R_i(t)$, using the same notation as that in Section 2 and assuming the censoring variables to be i.i.d.. The observed Fisher information at calendar time $t$ is

$$V_n(t) = \sum_{i=1}^{n} \delta_i(t) \left[ \sum_{j \in R_i(t)} z_j^2 / |R_i(t)| - \left\{ \sum_{j \in R_i(t)} z_j / |R_i(t)| \right\}^2 \right], \tag{11}$$

which provides an estimate of the null variance of $S_n(t)$. Suppose one uses a repeated significance test that rejects $H_0$ at the $j$th interim analysis $(1 \leq j \leq k)$ if

$$S_n(t_j)/V_n^{1/2}(t_j) \geq b_j \text{ or } S_n(t_j)/V_n^{1/2}(t_j) \leq a_j, \tag{12}$$

where $a_j < 0 < b_j$, and stops the trial at the first time $\tau \in \{t_1, \ldots, t_k\}$ when (12) occurs. Lai and Li [36] proposed to order the sample space of $(\tau, \Psi_\tau)$ by

$$(\tau_1, \Psi_{\tau_1}^{(1)}) \leq (\tau_2, \Psi_{\tau_2}^{(2)}) \quad \text{if and only if} \quad \Psi_{\tau_1 \wedge \tau_2}^{(1)} \leq \Psi_{\tau_1 \wedge \tau_2}^{(2)}, \tag{13}$$

where $\Psi_t = S(t)/V(t)$. They also proposed to estimate the unknown baseline distribution $G = 1 - e^{-\Lambda}$ by $\widehat{G} = 1 - e^{-\widehat{\Lambda}}$, where $\widehat{\Lambda}$ is Breslow's [6] estimator of the cumulative hazard function from all the data at the end of the trial:

$$\widehat{\Lambda}(s) = \sum_{i:Y_i(\tau) \leq s} \left\{ \delta_i(\tau) \middle/ \left( \sum_{j \in R_i(\tau)} e^{\widehat{\beta} z_j} \right) \right\}, \tag{14}$$

in which $\widehat{\beta}$ is Cox's [12] estimate of $\beta$ that maximizes the partial likelihood at time $\tau$. Since the $\xi_i$ are censored by $\min\{Y_i, (\tau - T_i)^+\}$, the distribution $C$ of $\xi_i$ can be estimated by the Kaplan-Meier estimator $\widehat{C}$. Let

$$\widehat{p}(\beta) = P\{(\tau^{(\beta)}, \Psi_{\tau^{(\beta)}}^{(\beta)}) > (\tau, \Psi_\tau)_{\mathrm{obs}}\}, \tag{15}$$

where the superscript $(\beta)$ means that the observations are generated by the proportional hazards model with baseline distribution $\widehat{G}$, censoring distribution $\widehat{C}$ and regression parameter $\beta$. As shown by Lai and Li [36], the confidence set

$$\{\beta : \alpha < \widehat{p}(\beta) < 1 - \alpha\} \tag{16}$$

has coverage probability $1 - 2\alpha + O(n^{-1/2})$ and is usually an interval under certain regularity conditions.

Lai and Li [36] used direct Monte Carlo to compute $\widehat{p}(\beta)$, and the confidence interval thus constructed is computationally intensive. The reason is that a large number of simulations are needed to compute $\widehat{p}(\beta)$ for each $\beta$, and a sample of survival times needs to be generated for each simulation. To reduce the computation time, we can use importance sampling by re-writing $\widehat{p}(\beta)$ as

$$\widehat{p}(\beta) = E_{\widehat{\beta}}\left[ \frac{L(\beta)}{L(\widehat{\beta})} I_{\{(\tau^{(\beta)}, \Psi_{\tau^{(\beta)}}^{(\beta)}) > (\tau, \Psi_\tau)_{\mathrm{obs}}\}} \right], \tag{17}$$

in which $L(\cdot)$ is the full likelihood at time $\tau$. This importance sampling technique provides a one-pass algorithm that only needs to generate data *once* under $P_{\widehat{\beta}}$ (after tilting $P_\beta$ to $P_{\widehat{\beta}}$), instead of having to generate data for each possible value $\beta$ of the confidence set in the direct Monte Carlo approach. For the Beta-Blocker Heart Attack Trial, Lai and Li [36] computed the confidence interval for the hazard ratio by direct Monte Carlo, which took over one day on a computer with Pentium 4 CPU 2.4GHz and 1024MB of RAM, in contrast to about 2.5 hours to compute the hybrid resampling confidence interval via (17) on the same computer.

An alternative approach that is commonly used in the literature is to use the space-time Brownian motion approximation of $(S_n(t), V_n(t))$ (see Jones and Whitehead [31] and Siegmund [59]) to which Siegmund's ordering can be applied. The following example, however, shows that the confidence intervals for $\beta$ constructed by this approach may have quite inaccurate coverage errors.

**Example 2.** Consider a time-sequential trial in which $n = 350$ subjects enter the trial uniformly during a 3-year recruitment period and are randomized to treatment or control with probability $\frac{1}{2}$. The trial is designed to last for a maximum of $t^* = 5.5$ years, with interim analyses after 1 year and every 6 months thereafter. The logrank statistic is used to test $H_0 : \beta = 0$ at each data monitoring time $t_j$ $(j = 1, \dots, 10)$ and the test is stopped at the smallest $t_j$ such that

$$V_n(t_j) \geq 55, \text{ or } V_n(t_j) \geq 11 \text{ and } |S_n(t_j)|/V_n^{\frac{1}{2}}(t_j) \geq 2.85, \tag{18}$$

or at $t_{10}(= t^*)$ when (18) does not occur, where $V_n(t)$ is defined by (11). If the test stops with $V_n(t_j) \geq 55$ or at $t^*$, reject $H_0$ if $|S_n(t^*)|/V_n^{\frac{1}{2}}(t^*) \geq 2.05$. Also reject $H_0$ if the second event in (18) occurs for some $j < 10$. The threshold 2.05 for the final analysis at $t^*$ is chosen so that the Type I error probability of the test is approximately 5% using the Brownian motion approximation; see Siegmund [59], page 132. When the space-time Brownian motion approximation of $(S_n(t), V_n(t))$ is applied in conjunction with Siegmund's ordering, thereby yielding score-based confidence intervals, the lower and upper coverage errors for nominal value $\alpha = 5\%$ are 4.45% and 5.05% for $\beta = 0$, 4.65% and 0.35% for $\beta = \log \frac{2}{3}$, and 5.75% and 0.00% for $\beta = \log \frac{1}{2}$ with 2000 simulations. This shows that the Brownian motion approximation does not provide an adequate approximation unless $\beta$ is very close to 0. The problem is that it fails to incorporate calendar time besides the information time. In contrast, Lai and Li's [36] hybrid resampling confidence set (16) has coverage errors 4.45% and 4.55% for $\beta = 0$, 5.25% and 5.35% for $\beta = \log \frac{2}{3}$, and 5.05% and 4.05% for $\beta = \log \frac{1}{2}$.